\magnification1200
\input amstex
\documentstyle{amsppt}

\def\ti{\times}
\def\IN{\Bbb N}
\def\IR{\Bbb R}
\def\St{\Cal St}
\def\V{\Cal V}
\def\W{\Cal W}
\def\w{\omega}
\def\C{\Cal C}
\def\diam{\operatorname{diam}}
\def\LRa{\Leftrightarrow}
\def\A{\Cal A}
\def\M{\Cal M}
\def\ol{\overline}
\def\a{\alpha}
\def\I{\Cal I}
\def\id{\operatorname{id}}
\def\s{\sigma}
\def\U{\Cal U}
\def\Ra{\Rightarrow}
\def\F{\Cal F}
\def\em{\it}
\def\ord{\operatorname{ord}}
\NoBlackBoxes
\hsize170truemm \vsize230truemm

\topmatter
\title
(Metrically) quarter-stratifiable spaces and their applications\\
in the theory of separately continuous functions
\endtitle
\author
Taras Banakh
\endauthor
\rightheadtext{(Metrically) quarter-stratifiable spaces and their
applications} \abstract{We introduce and study (metrically)
quarter-stratifiable spaces and then apply them to generalize
Rudin and Kuratowski-Montgomery theorems about the Baire and Borel
complexity of separately continuous functions. }
\endabstract
\address{Department of Mathematics, Lviv National University, Universytetska 1,
Lviv, 79000, Ukraina}
\endaddress
\email tbanakh\@franko.lviv.ua
\endemail
\endtopmatter

\document
{\baselineskip17pt
The starting point for writing this paper was
the desire to improve the results of V.K.~Mas\-lyu\-chenko et al.
\cite{MMMS}, \cite{MS}, \cite{KM}, \cite{KMM} who generalized a
classical theorem of W.Rudin \cite{Ru} which states that every
separately continuous function $f:X\ti Y\to \IR$ on the product of
a metrizable space $X$ and a topological space $Y$ belongs to the
first Baire class. It was proven in \cite{MMMS} that the
metrizability of $X$ in the Rudin theorem can be weakened to the
$\s$-metrizability and paracompactness of $X$. A subtle analysis
of Rudin's original proof reveals that this theorem is still valid
for a much wider class of spaces $X$. These spaces are of
independent interest, so we decided to give them a special name
--- metrically quarter-stratifiable spaces. (Metrically)
quarter-stratifiable spaces are introduced and studied in details
in the first three sections of this paper, where we investigate
relationships between the class of (metrically)
quarter-stratifiable spaces and other known classes of generalized
metric spaces. It turns out that each semi-stratifiable space is
quarter-stratifiable (this is a reason for the choice of the term
``quarter-stratifiable''), while each quarter-stratifiable
Hausdorff space has $G_\delta$-diagonal. Because of this, the
class of quarter-stratifiable spaces is ``orthogonal'' to the
class of compacta --- their intersection contains only metrizable
compacta. The class of quarter-stratifiable spaces is quite wide
and has many nice inheritance properties. Moreover, every
(submetrizable) space with $G_\delta$-diagonal is homeomorphic to
a closed subset of a (metrically) quarter-stratifiable
$T_1$-space. The following diagram describes the interplay between
the class of (metrically) quarter-stratifiable spaces and other
classes of generalized metric spaces in the framework of Hausdorff
spaces.
\newpage

}

{\parindent0pt
\hskip2.5cm{(metrizable)}

\hskip2.5cm{$\swarrow$ \hskip1cm $\searrow$}

\hskip25pt($\sigma$-metrizable)\hskip10pt  (stratifiable) \hskip15pt $\to$ \hskip15pt(paracompact quarter-stratifiable)

\hskip2.7cm $\searrow$ \hskip1cm $\downarrow$

\hskip3cm ($\sigma$-space)

\hskip4cm $\downarrow$

\hskip3.5cm(semi-stratifiable) \hskip100pt $\downarrow$

\hskip3.3cm$\swarrow$ \hskip.7cm --$\!\!\!\uparrow$ \hskip1.3cm   $\searrow$

\hskip10pt(hereditarily quarter-stratifiable) \hskip15pt (perfect)

\hskip3.8cm$\downarrow$\hskip3.3cm --$\!\!\!\!\nwarrow$

\hskip3cm(quarter-stratifiable) \hskip10pt $\gets$ \hskip10pt (metrically quarter-stratifiable)

\hskip2.8cm $\swarrow$\hskip2cm --$\!\!\!\uparrow$

\hskip1.5cm($G_\delta$-diagonal) $\longleftarrow$  (paracompact $G_\delta$-diagonal) \hskip30pt $\downarrow$

\hskip2cm$\uparrow$ \hskip4cm                     $\downarrow$

\hskip10pt (submetrizable) $\longleftrightarrow$ (closed subspace of a metrically quarter-stratifiable space)
}
\vskip15pt

In the last three sections we apply metrically quarter-stratifiable spaces for
generalizing classical theorems of Rudin \cite{Ru} and Kuratowski-Montgomery
\cite{Ku$_1$}, \cite{Ku$_2$}, \cite{Mo}.

\heading 1. Quarter-stratifiable spaces
\endheading

Before introducing quarter-stratifiable spaces we recall some
concepts from the theory of generalized metric spaces, see
\cite{Gr}. All topological spaces considered in this paper are
$T_1$-spaces, all maps (unlike to functions) are continuous.

\definition{1.1. Definition}
A topological space $X$ is \roster
\item {\em perfect} if each open set in $X$ is an $F_\sigma$-set;
\item {\em submetrizable} if $X$ admits a {\em condensation} (i.e.,
a bijective continuous map) onto a metrizable space;
\item {\em $\s$-metrizable} if $X$ is covered by a countable collection of
closed metrizable subspaces;
\item a {\em  $G_\delta$-diagonal} if the diagonal is a $G_\delta$-set in
the square $X\times X$;
\item a {\em  $\sigma$-space} if $X$ has a $\sigma$-discrete network;
\item {\em developable} if there exists a sequence $\{\U_n\}_{n\in\IN}$ of open
covers of $X$ such that for each $x\in X$ the set
$\{\St(x,\U_n)\}_{n\in\IN}$ is a base at $x$, where
$\St(x,\U_n)=\cup\{U\in\U_n: x\in U\}$;
\item a {\em  Moore space} if $X$ is regular and developable;
\item {\em semi-stratifiable} (resp. {\em stratifiable}) if there exists a
function $G$ which assigns to each $n\in\IN$ and a closed subset
$H\subset X$, an open set $G(n,H)$ containing $H$ such that
\itemitem{(i)} $H=\bigcap_{n\in\IN}G(n,U)$ (resp.
$H=\bigcap_{n\in\IN}\overline{G(n,U)}$)  and
\itemitem{(ii)} $G(n,K)\supset G(n,H)$ for every closed subset $K\supset
H$ and $n\in\IN$.
\endroster
\enddefinition

Semi-stratifiable spaces admit the following characterization, see
\cite{Gr, 5.8}.

\proclaim{1.2. Theorem} A topological space $(X,\tau)$ is
semi-stratifiable if and only if there exists a function $g:\IN\ti
X\to\tau$ such that for every $x\in X$
\roster
\item $\{x\}=\bigcap_{n\in\IN}g(n,x)$
\item $(x\in g(n,x_n),\; n\in\IN)\Ra (x_n\to x)$.
\endroster
\endproclaim

Weakening the first condition to $\bigcup_{x\in X}g(n,x)=X$ leads
to the definition of a quarter-stratifiable space --- the
principal concept in this paper.

\definition{1.3. Definition}
A topological space $(X,\tau)$ is called {\em
quarter-stratifiable} if there exists a function $g:\IN\ti
X\to\tau$ (called {\em quarter-stratifying function}) such that
\roster
\item $X=\bigcup_{x\in X}g(n,x)$ for every $n\in\IN$;
\item $(x\in g(n,x_n),\; n\in\IN)\Ra (x_n\to x)$.
\endroster
A space $X$ is {\em hereditarily quarter-stratifiable} if every
subspace of $X$ is quarter-stratifiable.
\enddefinition

Since the semi-stratifiability is a hereditary property, we get
that each semi-stratifiable space is hereditarily
quarter-stratifiable. As we shall see later, the converse is not
true: the Sorgenfrey arrow, being hereditarily
quarter-stratifiable, is not semi-stratifiable.

First, we give characterizations of the quarter-stratifiability in
terms of open covers as well as of lower semi-continuous
multivalued functions. We recall that a multivalued function
$\F:X\to Y$ is {\em lower semi-continuous} if for every open set
$U\subset Y$ its ``preimage'' $\F^{-1}(U)=\{x\in X:\F(x)\cap
U\ne\emptyset\}$ is open in $X$. For a topological space $X$ by
$X_d$ we denote $X$ endowed with the discrete topology.

\proclaim{1.4. Theorem} For a space $X$ the following statements
are equivalent: \roster
\item $X$ is quarter-stratifiable;
\item there exists a sequence $\{\U_n\}_{n\in\IN}$ of open covers of $X$
and a sequence $\{s_n:\U_n\to X\}_{n\in\IN}$ of functions such
that $s_n(U_n)\to x$ if $x\in U_n\in\U_n$, $n\in\IN$;
\item there exists a sequence $\{\F_n:X\to X_d\}_{n\in\IN}$ of lower
semi-continuous multivalued functions tending to the identity map
of $X$ in the sense that for every point $x\in X$ and every
neighborhood $O(x)\subset X$ of $x$ there is $n_0\in\IN$ such that
$\F_n(x)\subset O(x)$ for every $n\ge n_0$.
\endroster
\endproclaim

\demo{Proof} $(1)\Ra(2)$ Assume $(X,\tau)$ is quarter-stratifiable
and fix a function $g:\IN\ti X\to\tau$ from Definition 1.3. Let
$\U_n=\{g(n,x):x\in X\}$ and for every $U\in\U_n$ find $x\in X$
with $U=g(n,x)$ and let $s_n(U)=x$. Clearly that the so-defined
sequences $\{\U_n\}$ and $\{s_n\}$ satisfy our requirements.
\vskip5pt

$(2)\Ra(3)$ For every $n\in\IN$ define $\F_n:X\to X_d$ letting
$\F_n(x)=\{s_n(U):x\in U\in\U_n\}$ for $x\in X$. Clearly, $\F_n$
is lower-semicontinuous and $\F_n(x)\to x$ for each $x\in X$.
\vskip5pt

$(3)\Ra(1)$ For each $(n,x)\in\IN\ti X$ let $g(n,x)=\F_n^{-1}(x)$
and check that the so-defined function $g:\IN\times X\to\tau$
fits.\qed
\enddemo

Given a topological space $X$ by $l(X)$ and $d(X)$ we denote
respectively the Lindel\"of number and the density of $X$.

Like other classes of generalized metric spaces, the class of
quarter-stratifiable spaces has many nice inheritance properties.

\proclaim{1.5. Theorem} Let $X$ be a quarter-stratifiable space.
\roster
\item Every open subspace and every retract of $X$ is quarter-stratifiable.
\item If $X$ is Hausdorff, then it has a $G_\delta$-diagonal and countable
pseudo-character.
\item If $X$ is Hausdorff, then every paracompact \v Cech-complete subspace
and every countably compact subspace of $X$ is metrizable.
\item $d(X)\le l(X)$.
\item If $f:Y\to X$ is a condensation with sequentially continuous inverse
$f^{-1}$, then $Y$ is quarter-stratifiable.
\item If $f:X\to Y$ is a finite-to-one open surjective map, then the space
$Y$ is quarter-stratifiable.
\endroster
\endproclaim

\demo{Proof} Using Theorem~1.4, fix  a map $g:\IN\times X\to\tau$,
a sequence $\{\F_n:X\to X_d\}_{n\in\IN}$ of lower-semicontinuous
multivalued functions tending to the identity map of $X$, and
sequences $\{\U_n\}_{n\in\IN}$ of open covers and functions
$\{s_n\}_{n\in\IN}$ satisfying the conditions of Definition~1.3
and Theorem~1.4.
\vskip5pt

1) Given an open subspace $Y$ of $X$, let $\V_n=\{U\cap
Y:U\in\U_n\}$, $n\in\IN$. For each $V\in\V_n$ find a set
$U(V)\in\U_n$ with $V=U(V)\cap Y$ and let $\tilde
s_n(V)=s_n(U(V))$ if $s_n(U)\in Y$ and $\tilde s_n(V)$ be any
point of $Y$, otherwise. Clearly, $\V_n$ is an open cover of $Y$.
Now fix any point $y$ of $Y$ and take a sequence $V_n\in\V_n$ with
$y\in V_n$ for $n\in\IN$. Since $s_n(U(V_n))\to y$, we get
$s_n(U(V_n))\in Y$ for all $n$ beginning from some $n_0$. Then
$\tilde s_n(V_n)=s_n(U(V_n))$ for $n\ge n_0$ and thus $\tilde
s_n(V_n)\to y$, which proves that $Y$ is a quarter-stratifiable
space.

Now suppose that a subspace $Y$ of $X$ is a retract of $X$ and let
$r:X\to Y$ be the corresponding retraction. One can readily prove
that $\{r\circ{\F}_n|Y:Y\to Y_d\}$ is a sequence of
lower-semicontinuous multivalued functions tending to the identity
map of $Y$. \vskip5pt

2) Assume that the quarter-stratifiable space $X$ is Hausdorff.
Let $G_n=\bigcup_{U\in\U_n}U\times U\subset X\times X$. We claim
that $\bigcap_{n\in\IN}G_n$ coincides with the diagonal of
$X\times X$. To show this, take any distinct points $x,y\in X$ and
take disjoint neighborhood $O(x),O(y)$ of $x$ and $y$,
respectively. It follows from the choice of the sequences
$\{\U_n\}$ and $\{s_n\}$ that there is $n\in\IN$ such that
$s_m(U_m)\in O(x)$ and $s_m(V_m)\in O(y)$ for all $m\ge n$ and
$x\in U_m\in\U_m$, $y\in V_m\in\U_m$. We claim that $(x,y)\notin
G_n$. Assuming the converse, we would find a set $U_n\in\U_n$ with
$(x,y)\in U_n\times U_n$. Then $x,y\in U_n$ and $s_n(U_n)\in
O(x)\cap O(y)$, a contradiction. Thus $X$ has a
$G_\delta$-diagonal, which implies the countability of the
pseudo-character of $X$. \vskip5pt

3) If $X$ is Hausdorff and $Y\subset X$ is a paracompact \v
Cech-complete or countably compact subspace of $X$, then $Y$ has a
$G_\delta$-diagonal and by \cite{Bo, 8.2} or \cite{Gr, 2.14} is
metrizable. \vskip5pt

4) To see that $d(X)\le l(X)$ we may take the subset
$D=\{s_n(U):U\in\V_n,n\in\IN\}$, where $\V_n$ is a subcover of
$\U_n$ of size $|\V_n|\le l(X)$ for $n\in\IN$. Clearly, $D$ is a
dense subset in $X$ of size $|D|\le l(X)$. \vskip5pt

5) If $f:Y\to X$ is a condensation with sequentially continuous
inverse $f^{-1}$, we may consider the function $g_Y$ assigning to
each pair $(n,y)\in\IN\times Y$ the open set $f^{-1}(g(n,f(y)))$
of $Y$. One may easily show that the so-defined function $g_Y$
turns $Y$ into a quarter-stratifiable space. \vskip5pt

6) Suppose $f:X\to Y$ is an open finite-to-one surjective map. It
can be shown that for every $n\in\IN$ the multivalued map
$f\circ\F_n\circ f^{-1}:Y\to Y_d$ defined by $f\circ\F_n\circ
f^{-1}=\bigcup_{x\in f^{-1}(y)}f(\F_n(x))$ for $y\in Y$ is
lower-semicontinuous. Moreover, the sequence $\{f\circ\F_n\circ
f^{-1}\}_{n\in\IN}$ tends to the identity map of $Y$. Then
according to Theorem~1.4,
the space $Y$ is quarter-stratifiable. \qed
\enddemo

Next, we show that the quarter-stratifiability is stable with
respect to certain set-theoretic operations.

\proclaim{1.6. Theorem}
\roster
\item The countable product of quarter-stratifiable spaces is
quarter-stratifiable.
\item A space $X$ is quarter-stratifiable, provided $X=A\cup B$ is a union
of two quarter-stratifiable subspaces of $X$, one of which is a
closed $G_\delta$-set in $X$.
\item A space $X$ is quarter-stratifiable, provided $X$ can be covered by
a $\sigma$-locally finite collection of quarter-stratifiable
closed $G_\delta$-subspaces of $X$.
\endroster
\endproclaim

\demo{Proof}

1) Let $X_i$, $i\in\IN$, be a sequence of quarter-stratifiable
spaces. For every $i\in\IN$, fix a point $*_i\in X_i$ and a
sequence $\{\F_{i,n}:X_i\to (X_i)_d\}_{n\in\IN}$ of
lower-semicontinuous multivalued functions tending to the identity
map of $X_i$ as $n\to\infty$. Let $X=\prod_{i\in\IN}X_i$ and for
every $n\in\IN$ identify the product $\prod_{i=1}^nX_i$ with the
subspace $\{(x_i)_{i\in\IN}:x_i=*_i$ for $i>n\}$ of $X$. For every
$n\in\IN$ define a multivalued function $\F_n:X\to X_d$ as
follows: $\F_n((x_i)_{i\in\IN})=\prod_{i=1}^n\F_{i,n}(x_i)$. It is
easy to see that each function $\F_n$ is lower semicontinuous and
the sequence $\{\F_n\}$ tends to $\id_X$. \vskip5pt

2) Suppose $X=A\cup B$, where $A,B$ are quarter-stratifiable
subspaces and $B$ is a closed $G_\delta$-set in $X$. Fix a
decreasing sequence $\{O_n(B)\}_{n\in\IN}$ of open subsets of $X$
with $B=\bigcap_{n\in\IN}O_n(B)$. By Theorem~1.5(1), the open
subspace $A\setminus B=X\setminus B$ of $A$ is
quarter-stratifiable and thus admits a function $g_{A\setminus
B}:\IN\times (A\setminus B)\to\tau$ into the topology $\tau$ of
$X$ such that $\bigcup_{x\in A\setminus B}g_{A\setminus
B}(n,x)\supset A\setminus B$ for each $n\in\IN$, and $x\in
g_{A\setminus B}(n,x_n)\;\Ra\; x_n\to x$ for each $x,x_n\in
A\setminus B$, $n\in\IN$.

Using the quarter-stratifiability of $B$ find a function
$g_B:\IN\times B\to\tau$ such that $B\subset\bigcup_{x\in
B}g_B(n,x)\subset O_n(B)$ for every $n\in\IN$, and $x\in
g_B(n,x_n)\;\Ra\;x_n\to x$ for each $x,x_n\in B$, $n\in\IN$. It is
easy to verify that the function $g:\IN\ti X\to\tau$ defined by
$$g(n,x)=\cases g_{A\setminus B}(n,x) &\text{if $x\in A\setminus
B$,}\\ g_B(n,x)& \text{if $x\in B$}
\endcases
$$ turns $X$ into a quarter-stratifiable space. \vskip5pt

3) Let $\{X_i\}_{i\in\I}$ be a $\sigma$-locally finite collection
of quarter-stratifiable closed $G_\delta$-subspaces of a space
$X$. Write $\I=\bigcup_{k\in\IN}\I_k$ so that the collection
$\{X_i\}_{i\in\I_k}$ is locally finite for every $k\in\IN$.
Without loss of generality, $\I_k\subset \I_{k+1}$ for every $k$.
For every $k$ fix an open cover $\W_k$ of $X$ whose any element
$W\in\W_k$ meets only finitely many of the sets $X_i$'s with
$i\in\I_k$. We may additionally require that the cover $\W_k$ is
inscribed into the cover $\W_{k-1}$ for $k>1$.

Since $X_i$'s are $G_\delta$-sets in $X$, for every $i\in\I$ we
may find a decreasing sequence $\{O_n(X_i)\}_{n\in\IN}$ of open
subsets $O_n(X_i)\subset\St(X_i,\W_n)=\{W\in\W_n:W\cap
X_i\notin\emptyset\}$ with $\bigcap_{n\in\IN}O_n(X_i)=X_i$.

Let $\le$ be any well-ordering of the index set $\I$ such that
$i<j$ for every $i,j\in I$ with $i\in I_k\not\ni j$ for some $k$.
It follows from the local finiteness of the collections
$\{X_i\}_{i\in\I_k}$ that for every $i\in\I$ the set
$\bigcup_{j<i}X_j$ is closed in $X$. Using the
quarter-stratifiability of open subsets of $X_i$'s, for every
$i\in\IN$ fix a sequence $\{\U_{i,n}\}_{n\in\IN}$ of open covers
of the set $Y_i=X_i\setminus\bigcup_{j<i}X_j$ and a sequence
$\{s_{i,n}:\U_{i,n}\to Y_i\}$ of functions such that
$s_{i,n}(U_n)\to y$ for every $Y_i\ni y\in U_n\in\U_{i,n}$,
$n\in\IN$. Without loss of generality, we may assume that each
$U\in\U_{i,n}$ is an open subset of $X$ lying in
$O_n(X_i)\setminus \bigcup_{j<i}X_j$. Now for every $n\in\IN$
consider the open cover $\U_n=\bigcup_{i\in\I}\U_{i,n}$ and the
function $s_n=\bigcup_{i\in\I}s_{i,n}:\U_n\to X$.

To show that the space $X$ is quarter-stratifiable, it suffices to
verify that for every $x\in X$ and a sequence $U_n\in\U_n$,
$n\in\IN$, with $x\in U_n$ we have $s_n(U_n)\to x$. Let
$i=\min\{j\in\I:x\in X_i\}$. Then $x\in Y_i$. For every $n\in\IN$
find $i_n\in\I$ with $U_n\in\U_{i_n,n}$. Since $x\in U_n\subset
O_{i_n}(X_i)\setminus\bigcup_{j<i_n}X_j$, we get $i_n\le i$ for
all $n\in\IN$. It follows from the choice of the cover $\W_n$ and
the sets $O_n(X_j)$ that
$\bigcap_{n=1}(\bigcup_{j<i}O_n(X_j))=\bigcup_{j<i}X_j$.
Consequently, $x\notin \bigcap_{n=1}^\infty
(\bigcup_{j<i}O_n(X_j))$ and there is $n_0\in\IN$ such that
$x\notin \bigcup_{j<i}O_n(X_j)$ for every $n\ge n_0$. Since $x\in
U_n\subset O_n(X_{i_n})$ for every $n\in\IN$, we conclude that
$i_n=i$ for $n\ge n_0$. Then $x\in U_n\in\U_{i,n}$ for all $n\ge
n_0$ and $s_n(U_n)=s_{i,n}(U_n)\to x$, which completes the proof
of the quarter-stratifiability of $X$.\qed
\enddemo

\definition{1.7. Remark} The $G_\delta$-condition in Theorem~1.6(2) is essential:
the one-point compactification $\alpha\Gamma$ of an uncountable
discrete space $\Gamma$ may be written as the union
$\alpha\Gamma=\{\infty\}\cup\Gamma$ of two metrizable subspaces,
one of which is closed. But $\alpha\Gamma$, being a non-metrizable
compactum, is not quarter-stratifiable, see Theorem~1.5(3). Yet,
we do not know the answer to the following question.
\enddefinition

\proclaim{1.8. Question} Is a (regular) space $X$
quarter-stratifiable if it is a union of two closed
quarter-stratifiable subspaces?
\endproclaim

According to Theorem~1.5(1), every open subspace of a
quarter-stratifiable space is quarter-stratifiable. It is not true
for closed subspaces of quarter-stratifiable spaces: every
$G_\delta$-diagonal space is homeomorphic to a closed subset of a
quarter-stratifiable $T_1$-space.

\proclaim{1.9. Theorem} Every space $X$ with $G_\delta$-diagonal
is homeomorphic to a closed subspace of a quarter-stratifiable
$T_1$-space $Y$ with $|Y\setminus X|\le\max\{l(X),d(X)\}$.
\endproclaim

\demo{Proof} Suppose $X$ is a space with $G_\delta$-diagonal. By
Theorem 2.2 of \cite{Gr}, $X$ admits a sequence
$\{\V_n\}_{n\in\IN}$ of covers such that
$\{x\}=\bigcap_{n\in\IN}\St(x,\V_n)$ for every $x\in X$. Replacing
each $\V_n$ by a suitable subcover, we may assume that $|\V_n|\le
l(X)$ for every $n$, and each cover $\V_n$ consists of non-empty
subsets of $X$. Let $D\subset X$ be a dense subset of $X$ with
$|D|\le d(X)$. For every $n\in\IN$ and every $V\in\V_n$ pick a
point $c(V)\in V\cap D$. Let $S_n=\{0\}\cup\{1/i:i\ge n\}$ denote
a tail of the convergent sequence $S_1$.

Consider the subset $Y=X\times\{0\}\cup\{(c(V),1/n):V\in\V_n,\;
n\in\IN\}\subset X\times S_1$. It is clear that $|Y\setminus
X|\le\min\{|\cup_{n\in\IN}\V_n|,|D|\}\le\min\{l(X),d(X)\}$.
Identify $X$ with the subset $X\times \{0\}$ of $Y$. Define a
topology $\tau$ on $Y$ letting $U\subset Y$ be open if and only if
$U\cap X$ is open in $X$ and for every $x\in U\cap X$ there is
$n_0\in\IN$ such that $\{(c(V),1/n):x\in V\in\V_n,\; n\ge
n_0\}\subset U$. Thus $X$ is homeomorphic to a closed subspace of
$Y$, while all points of $Y\setminus X$ are isolated. Since each
one-point subset of $Y$ is closed, $Y$ is a $T_1$-space. To show
that $Y$ is a quarter-stratifiable space, consider open covers
$\U_n=\{\{y\},(V\times S_n)\cap Y:y\in Y\setminus X,\;
V\in\V_n\}$, $n\in\IN$, of $Y$. For every $U\in \U_n$ let $$
s_n(U)=\cases y & \text{if $U=\{y\}$ for some $y\in Y\setminus
X$};\\ (c(V),1/n) & \text{if $U=(V\times S_n)\cap Y$ for
$V\in\V_n$}.
\endcases
$$ It follows from the choice of the topology $\tau$ on $Y$ that
$x\in U_n\in\U_n$, $n\in\IN$, implies $s_n(U_n)\to x$ \, i.e., $Y$
is a quarter-stratifiable $T_1$-space. \qed
\enddemo

\proclaim{1.10. Problem} Describe the class of subspaces of
regular (Tychonoff) quarter-stratifiable $T_1$-spaces.
\endproclaim

\heading 2. Metrically quarter-stratifiable spaces
\endheading

In this section we introduce and study metrically
quarter-stratifiable spaces forming a class, intermediate between
the class of paracompact quarter-stratifiable Hausdorff spaces and
the class of submetrizable quarter-stratifiable spaces. We start
with defining a quarter-stratifying topology.

\definition{2.1. Definition}
A topology $\tau'$ on a topological space $(X,\tau)$ is called
{\em quarter-stratifying} if there exists a quarter-stratifying
function $g:\IN\ti X\to\tau$ for $X$ such that $g(\IN\ti
X)\subset\tau'$.

A topological space $X$ is defined to be {\it metrically
quarter-stratifiable} if it admits a weaker metrizable
quarter-stratifying topology.
\enddefinition

The following theorem characterizes metrically
quarter-stratifiable spaces.

\proclaim{2.2. Theorem} For a space $X$ the following statements
are equivalent: \roster
\item $X$ admits a weaker metrizable quarter-stratifying topology;
\item $X$ admits a weaker paracompact Hausdorff quarter-stratifying topology;
\item there exists a weaker metrizable topology $\tau_m$ on $X$,
a sequence $\{\U_n\}_{n\in\IN}$ of covers of $X$ by $\tau_m$-open
subsets and a sequence $\{s_n:\U_n\to X\}_{n\in\IN}$ of functions
such that $s_n(U_n)\to x$ in $X$ if $x\in U_n\in\U_n$, $n\in\IN$;
\item there exists a weaker paracompact Hausdorff topology $\tau_p$ on $X$,
a sequence $\{\V_n\}_{n\in\IN}$ of covers of $X$ by $\tau_p$-open
subsets and a sequence $\{\tilde s_n:\V_n\to X\}_{n\in\IN}$ of
functions such that $\tilde s_n(V_n)\to x$ in $X$ if $x\in
V_n\in\V_n$, $n\in\IN$;
\endroster
\endproclaim

\demo{Proof} The equivalence $(1)\LRa(3)$ and  $(2)\LRa(4)$ can be
proved by analogy with the corresponding equivalences in Theorem
1.4; the implication $(3)\Ra (4)$ is trivial.

$(4)\Ra(3)$ Assume that $\tau_p$ is a weaker paracompact Hausdorff
topology on $X$, $\{U_n\}_{n\in\IN}$ is a sequence of
$\tau_p$-open covers of $X$, and $\{s_n:\U_n\to X\}_{n\in\IN}$ is
a sequence of functions such that $s_n(U_n)\to x$ in $X$ if $x\in
U_n\in\U_n$, $n\in\IN$. Let $l_1(\U_n)$ denote the Banach space of
all absolutely summable functions $f:\U_n\to\IR$ equipped with the
norm $\|f\|=\sum_{U\in\U_n}|f(U)|$.

Using the paracompactness of the topology $\tau_p\subset\tau$, for
every $n\in\IN$ find a partition of unity
$\{\lambda_U:X\to[0,1]\}_{U\in\U_n}$ such that
$\lambda_U^{-1}(0,1]\subset U$ for $U\in\U_n$. Observe that this
partition of unity can be seen as a continuous map $\Lambda_n:X\to
l_1(\U_n)$ acting as $\Lambda_n(x)=(\lambda_U(x))_{U\in\U_n}$ (the
continuity of $\Lambda_n$ follows from the local finiteness of the
cover $\{\lambda_U^{-1}(0,1]:U\in\U_n\}$~).

Let $\Lambda:X\to \prod_{n\in\IN}l_1(\U_n)$ denote the diagonal
product of the maps $\Lambda_n$. By analogy with the proof of
Theorem 1.5(2), show that the map $\Lambda$ is injective. Then the
weakest topology $\tau_m$ on $X$ for which the map $\Lambda$ is
continuous is metrizable. Observe that for each $n\in\IN$ and
$U\in\U_n$ the set $V(U)=\lambda_U^{-1}(0,1]$ is $\tau_m$-open.
Let $\V_n=\{V(U):U\in\U_n\}$. For every $V\in\V_n$ find $U\in\U_n$
with $V=V(U)$ and let $\tilde s_n(V)=s_n(U)$. One can easily check
that the metrizable topology $\tau_m$ on $X$ and the sequences
$\{V_n\}_{n\in\IN}$, $\{\tilde s_n\}_{n\in\IN}$ satisfy the
condition (3). \qed
\enddemo

In the following theorem we collect some elementary properties of
metrically quarter-stratifiable spaces.

\proclaim{2.3. Theorem} \roster
\item Each paracompact Hausdorff quarter-stratifiable space is metrically
quarter-stratifiable.
\item Each metrically quarter-stratifiable space is
submetrizable and quarter-stratifiable.
\item Every open subspace of a metrically quarter-stratifiable
space is metrically quarter-stratifiable.
\item If $X$ is a metrically quarter-stratifiable space and
$f:Y\to X$ is a condensation with sequentially continuous inverse
$f^{-1}$, then the space $Y$ is metrically quarter-stratifiable.
\item The product $X=\prod_{n\in\IN}X_n$ of a countable collection of
metrically quarter-stratifiable spaces $X_n$, $n\in\IN$, is a
metrically quarter-stratifiable space.
\item A space $X$ is metrically quarter-stratifiable if there exist a weaker
paracompact Hausdorff topology $\tau_p$ and two subspaces
$A,B\subset X$ such that $A\cup B=X$, $A$ is a closed
$G_\delta$-set in $(X,\tau_p)$ and the topologies induced by the
topology $\tau_p$ on $A$ and $B$ are quarter-stratifying.
\item A space $X$ is metrically quarter-stratifiable if there exist
a weaker paracompact Hausdorff topology $\tau_p$ on $X$ and a
$\sigma$-locally finite cover $\C$ of  $(X,\tau_p)$ by closed
$G_\delta$-subspaces such that the topology induced by $\tau_p$ on
each $C\in \C$ is quarter-stratifying (with respect to the
original topology of $C$);
\item A space $X$ is metrically quarter-stratifiable if there is a
weaker metrizable topology $\tau_m$ on $X$ and a cover $\C$ of
$X$, well-ordered by the inclusion relation, such that the
topology $\tau_m$ induces the original topology on each $C\in\C$;
\item Every submetrizable space is homeomorphic to a metrically quarter-
stratifiable space.
\endroster
\endproclaim

\demo{Proof} All statements (except for (8)) easily follow from
definitions
 or can be proved by analogy with the
corresponding properties of quarter-stratifiable spaces.

To prove the eighth statement, fix a continuous metric $d$ on $X$
and a well-ordered (by the inclusion relation) cover $\C$ of $X$
such that $d$ induces the original topology on each $C\in\C$. Let
$\U_n$ denotes the collection of all open $1/n$-balls with respect
to the metric $d$. For every $U\in\U_n$ let $C(U)$ be the smallest
set $C\in\C$ meeting $U$ and let $s_n(U)\in U\cap C(U)$ be any
point.

Fix any $x\in X$ and a sequence $x\in U_n\in\U_n$, $n\in\IN$. We
have to show that $s_n(U_n)\to x$. Let $C(x)$ be the smallest set
$C\in\C$ containing the point $x$. Since $x\in C(x)\cap U_n$, we
conclude that $s_n(U_n)\in C(x)$. By the choice of the sets $U_n$,
$\diam(U_n)\le 2/n$. Then $d(x,s_n(U_n))\le \diam(U_n)\le 2/n$,
$n\in\IN$, and thus $s_n(U_n)\to x$ (because $d$ induces the
topology of $C(x)$). \qed
\enddemo

According to Theorem 2.3(2) each metrically quarter-stratifiable
space is submetrizable and quarter-stratifiable. We do not know if
the converse is also true.

\proclaim{2.4. Question} Is every submetrizable
quarter-stratifiable space metrically quarter-strati\-fi\-able?
\endproclaim

\heading 3. Some Examples
\endheading

In this section we collect some examples exposing the difference
between the class of (metrizable) quarter-stratifiable spaces and
other classes of generalized metric spaces. It is known that every
Moore space is semi-stratifiable and each collectively normal
Moore space is metrizable, see \cite{Gr}. Yet, there exist Moore
spaces which are not submetrizable, see \cite{Ny}. Thus we have

\proclaim{3.1. Example} There exists a quarter-stratifiable
(Moore) space which is not metrically quarter-stratifiable.
\endproclaim

Theorems~1.9 and 2.3(9) show that the class of
quarter-stratifiable spaces is much wider than the class of
semi-stratifiable spaces. Using these theorems many wild
(metrically) quarter-stratifiable spaces may be constructed. But
in general, so-constructed spaces are not hereditarily
quarter-stratifiable.

\proclaim{3.2. Example} The Sorgenfrey arrow $Z$ is not
semi-stratifiable but every subspace of $Z$ is Lindel\"of,
separable, and metrically quarter-stratifiable.
\endproclaim

\demo{Proof} Recall that the Sorgenfrey arrow $Z$ is the
semi-interval $[0,1)$ endowed with the topology generated by the
base consisting of all semi-intervals $[a,b)$ where $0\le a<b\le
1$, see \cite{En$_1$, 1.2.2}. Since $Z$ embeds into a linearly
ordered space (``two arrows'' of Aleksandrov \cite{En$_1$,
3.10.C}), $Z$ is monotonically normal, see \cite{Gr, 5.21}. If $Z$
would be semi-stratifiable, then $Z$, being monotonically normal,
would be stratifiable according to Theorem 5.16 of \cite{Gr}.
Since each stratifiable space is a $\s$-space \cite{Gr, 5.9} and
Lindel\"of $\s$-spaces have countable network \cite{Gr, 4.4}, this
would imply that $Z$ has countable network, which is a
contradiction. Therefore, the Sorgenfrey arrow $Z$ is not
semi-stratifiable.

Next, we show that every subspace $X$ of $Z$ is
quarter-stratifiable. For every $n\in\IN$ consider the finite open
cover $\U_n=\{X\cap[\frac{k-1}n,\frac kn):0< k\le n\}$ of $X$. For
every element $U\in\U_n$ choose a point $s_n(U)\in X$ as follows.
If the set $\uparrow U=\{x\in X:x\ge\sup U\}$ is not empty, then
let $s_n(U)$ be any point of $X$ with $s_n(U)<\inf \uparrow
U+\frac1n$; otherwise, let $s_n(U)$ be any point of $X$. It can be
shows that $s_n(U_n)\to x$ for every $x\in X$ and every sequence
$x\in U_n\in\U_n$, $n\in\IN$. Thus the space $X$ is
quarter-stratifiable. It is well known that the Sorgenfrey arrow
is hereditarily Lindel\"of and hereditarily separable which
implies that the space $X$ is Linde\"of and separable. Then $X$ is
also paracompact (\cite{En$_1$, 5.1.2}) and hence, being
quarter-stratifiable, is metrically quarter-stratifiable. \qed
\enddemo

Since the metrical quarter-stratifiability is productive, we get

\proclaim{3.3. Example} The square of the Sorgenfrey arrow is
metrically quarter-stratifiable but not normal.
\endproclaim

Another example of a metrically quarter-stratifiable non-normal
space is the Nemytski plane. Unlike to (semi-)stratifiable spaces
metrically quarter-stratifiable spaces need not be perfect (or
paracompact).


\proclaim{3.4. Example} A separable zero-dimensional metrically
quarter-stratifiable Tychonov space which is neither perfect nor
Lindel\"of.
\endproclaim

\demo{Proof} Consider the Cantor cube $2^\w=\{0,1\}^\w$, i.e., the
set of all functions $\w\to\{0,1\}$. There is a natural partial
order on $2^\w$: $x\le y$ iff $x(i)\le y(i)$ for each $i\in\w$.
For $k\in\{0,1\}$ let $2^\w_k=\{(x_i)_{i\in\w}: \exists n\in\w$
with $x(i)=k$ for all $i\ge n\}$. For a point $x\in 2^\w$ and
$n\in\w$ let $U(x,n)=\{y\in 2^\w: y\ge x$ and $y(i)=x(i)$ for
$i\le n\}$. On the set $2^\w$ consider the topology $\tau$
generated by the base $\{U(x,n):x\in 2^\w,\;n\in\w\}$. Clearly,
the space $X=(2^\w,\tau)$ is Tychonoff, zero-dimensional, and
separable ($2^\w_1$ is a countable dense set in $X$).

Let us show that the space $X$ is metrically quarter-stratifiable.
For every $n\in\w$ identify $2^n$ with the subset $\{(x_i)\in
2^\w: x_i=0$ for $i\ge n\}\subset 2^\w$ and consider the finite
cover $\U_n=\{U(x,n):x\in 2^n\}$ of $2^\w$. It is clear that each
set $U(x,n)\in\U_n$ is open is the product topology of $2^\w$
which is metrizable. To each $U\in\U_n$ assign the point
$s_n(U)=\max U\in 2^\w_1$. One may easily verify that $s_n(U_n)\to
x$ for every $x\in X$ and a sequence $x\in U_n\in\U_n$, $n\in\w$.
Thus the space $X$ is metrically quarter-stratifiable.

Now we show that the space $X$ is not perfect. Note that $2^\w_0$
is a closed subset of $X$. We claim that it is not a
$G_\delta$-subset of $X$. Assume conversely that
$2^\w_0=\bigcap_{n\in\IN}O_n$, where $O_n$ are open subsets of
$X$. Since the topology of $X$ coincides with the usual product
topology at every point of the set $2^\w_0$, we may assume that
each set $O_n$ is open in the product topology of $2^\w$. Then by
the Baire theorem, the intersection $\bigcap_{n\in\IN}O_n$ is a
dense $G_\delta$-subset in $2^\w$, which contains a point $x\notin
2^\w_0$, a contradiction.

By Theorem~2.3(3) the open subspace $Y=X\setminus 2^\w_0$ of $X$
is metrically quarter-stratifiable. Since $X$ is regular and
$2^\w_0$ is not a $G_\delta$-set in $X$, we conclude that the
space $Y$ is not Lindel\"of.

Then the product $X\times Y$ is a separable Tychonov
zero-dimensional metrically quarter-stratifiable space which is
neither perfect nor Lindel\"of. \qed
\enddemo

\proclaim{3.5. Question} Is there a hereditarily
quarter-stratifiable space which is not perfect?
\endproclaim

We recall that a topological space $X$ is {\it finally compact} if
$l(X)\le\aleph_0$, i.e., every open cover of $X$ admits a
countable subcover. Finally compact regular spaces are called {\em
Lindel\"of}.

\proclaim{3.6. Example} A hereditarily finally compact
submetrizable uncountable space whose uncountable subspaces are
neither separable nor quarter-stratifiable.
\endproclaim

\demo{Proof} Let $X$ be the interval $(0,1)$ with the topology
generated by the sets of the form $(a,b)\setminus C$, where $0\le
a<b\le1$ and $C$ is a countable set. Clearly, the space $X$ is
submetrizable and hereditarily finally compact. Yet, every
countable subset of $X$ is closed which implies that every
uncountable subspace $Y$ of $X$ is not separable. Since finally
compact quarter-stratifiable spaces are separable, see
Theorem~1.6(4), we conclude that the space $Y$ is not
quarter-stratifiable.\qed
\enddemo

The space constructed in Example~3.6 is not regular. There is also
a regular submetrizable non-quarter-stratifiable space.

Given a subset $B\subset\IR$ denote by $\IR_B$ the real line $\IR$
endowed with the topology consisting of the sets $U\cup A$, where
$U$ open in $\IR$ and $A\subset B$, see \cite{En$_1$, 5.1.22 and
5.5.2}. It is known that for every $B\subset \IR$ the space
$\IR_B$ is hereditarily paracompact; $\IR_B$ is perfect if and
only if $B$ is a $G_\delta$-set in $\IR$. If $B\subset\IR$ is a
Bernstein set, then the space $\IR_B$ is Lindel\"of, see
\cite{En$_1$, 5.5.4}. Recall that a subset $B\subset\IR$ is called
a {\it Bernstein set} if $C\cap B\ne\emptyset\ne C\setminus B$ for
every uncountable compactum $C\subset\IR$. Bernstein sets can be
easily constructed by transfinite induction, see \cite{En$_1$,
5.4.4}.

\proclaim{3.7. Example} If $B\subset\IR$ is not $\s$-compact, then
the space $\IR_B$ is not quarter-stratifiable. Consequently, if
$B\subset \IR$ is the set of irrationals, then $\IR_B$ is a
hereditarily paracompact submetrizable perfect zero-dimensional
space which is neither separable nor quarter-stratifiable; if
$B\subset\IR$ is a Bernstein set, then $\IR_B$ is a submetrizable
hereditarily paracompact Lindel\"of zero-dimensional space which
is neither perfect nor separable nor quarter-stratifiable.
\endproclaim

\demo{Proof} Assume that $B\subset\IR$ is not $\s$-compact.
Suppose the space $\IR_B$ is quarter-stratifiable and fix a
sequence $\{\U_n\}_{n\in\IN}$ of open covers of $\IR_B$ and a
sequence $\{s_n:\U_n\to \IR_B\}_{n\in\IN}$ of functions such that
$s_n(U_n)\to x$ for every $x\in U_n\in\U_n$, $n\in\IN$. Since the
real line is hereditarily Lindel\"of and the topology of $\IR_B$
coincides with the usual topology at points of the set
$\IR\setminus B$, we may find countable subcovers
$\V_n\subset\U_n$, $n\in\IN$, of the set $\IR\setminus B$. Observe
that $C=\{s_n(V):V\in\V_n,\;n\in\IN\}$ is countable while the
intersection $G=\bigcap_{n\in\IN}\St(\IR\setminus B,\V_n)$ is a
$G_\delta$-set in the usual topology of $\IR$. Since $B$ is not a
$\sigma$-compact, there is a point $x\in B\cap G\setminus C$.
Choose for every $n\in\IN$ an element $V_n\in\V_n$ with $V_n\ni
x$. Then $x\ne s_n(V_n)$ for all $n$ and by the definition of the
topology of $\IR_B$, $s_n(V_n)\not\to x$. This contradiction shows
that the space $\IR_B$ is not quarter-stratifiable.\qed
\enddemo

\heading 4. A generalization of Kuratowski-Montgomery Theorem
\endheading

In this section we apply metrically quarter-stratifiable spaces to
generalize a theorem of Kuratowski \cite{Ku$_2$} and Montgomery
\cite{Mo} which asserts that for metrizable spaces $X,Y,Z$ a
function $f:X\ti Y\to Z$ is Borel measurable of class $\a+1$, $\a$
a countable ordinal, if $f$ is continuous with respect to the
first variable and Borel measurable of class $\a$ with respect to
the second variable.

At first, we recall the definitions of the multiplicative and
additive Borel classes. Given a topological space $X$ let
$\A_0(X)$ and $\M_0(X)$ denote the classes of all open and all
closed subsets of $X$, respectively. Assuming that for a countable
ordinal $\a$ the classes $\A_\beta(X)$, $\M_\beta(X)$, $\beta<\a$,
are already defined, let $\A_a(X)=\{\bigcup_{n=1}^\infty M_n:
M_n\in\bigcup_{\beta<\a}\M_\beta(X)$ for all $n\in\IN\}$ and
$\M_\a(X)=\{M\subset X: X\setminus M\in\A_\a(X)\}$.

We say that a function $f:X\to Y$ between topological spaces is
Borel measurable of class $\a$ if the preimage $f^{-1}(U)$ of any
open set $U\subset Y$ belongs to the additive class $\A_\a(X)$ (of
course, this is equivalent to saying that the preimage $f^{-1}(F)$
of any closed subset $F\subset Y$ belongs to the multiplicative
class $\M_\a(X)$).

The set of all Borel measurable functions $f:X\to Y$ of class $\a$
is denoted by $H_\a(X,Y)$. For topological spaces $X,Y,Z$ and a
countable ordinal $\a$ let $CH_\a(X\ti Y,Z)$ denote the set of all
functions $f:X\ti Y\to Z$ which are continuous with respect to the
first variable and are Borel measurable of class $\a$ with respect
to the second variable. The Kuratowski-Montgomery Theorem states
that $CH_\a(X\ti Y,Z)\subset H_{\a+1}(X\ti Y,Z)$ for metrizable
spaces $X,Y,Z$ and a countable ordinal $\a$.

A subset $A$ of a topological space $X$ is called a {\it
$\ol{G}_\delta$-set} if
$A=\bigcap_{n\in\IN}U_n=\bigcap_{n\in\IN}\bar U_n$ for some open
sets $U_n\subset X$, $n\in\IN$. Observe that every
$\ol{G}_\delta$-set is a closed $G_\delta$-set and every closed
$G_\delta$-set in a normal space is a $\ol{G}_\delta$-set.
Consequently, every closed subset of a perfectly normal space is a
$\ol{G}_\delta$-set.

\proclaim{4.1. Theorem} The inclusion $CH_\a(X\ti Y,Z)\subset
H_{\a+1}(X\ti Y,Z)$ holds for any countable ordinal $\a$, any
metrically quarter-stratifiable space $X$, any topological space
$Y$, and any topological space $Z$ whose every closed subset is a
$\ol{G}_\delta$-set.
\endproclaim

\demo{Proof} Using Theorem 2.2, fix a weaker metrizable topology
$\tau_m$ on the metrically quarter-stratifiable space $X$, a
sequence $\{\U_n\}_{n\in\IN}$ of $\tau_m$-open covers of $X$ and a
sequence $\{s_n:\U_n\to X\}_{n\in\IN}$ of functions such that
$s_n(U_n)\to x$ if $x\in U_n\in\U_n$, $n\in\IN$. According to the
classical Stone Theorem \cite{St}, we may assume that each cover
$\U_n$ is locally finite and $\s$-discrete.

Fix any function $f\in CH_\a(X\ti Y,Z)$. To show that $f\in
H_{\a+1}(X\ti Y,Z)$ we have to verify that the preimage
$f^{-1}(F)$ of any closed subset $F\subset Z$ belongs to the class
$\M_{\a+1}(X\ti Y)$.

The set $F$, being $\ol{G}_\delta$-set, can be written as
$F=\bigcap_{n=1}^\infty \ol{W}_n$, where $W_n\subset W_{n-1}$ are
open neighborhoods of $F$. For every $n\in\IN$ and $U\in\U_n$ let
$F_U:Y\to Z$ be the map defined by $f_U(y)=f(s_n(U),y)$ for $y\in
Y$. It is clear that $f_U\in H_\a(Y,Z)$ which yields
$f^{-1}_U(W_m)\in\A_\a(Y)$ for every $m$. Since the set $U$ is
functionally open in $X$, we get $U\ti f_U^{-1}(W_m)\in\A_\a(X\ti
Y)$, $m\in\IN$. Using this observation and the $\s$-discreteness
of the cover $\U_n$, show that
$A_{n,m}=\bigcup\limits_{U\in\U_n}U\ti f_U^{-1}(W_n)\in\A_\a(X\ti
Y)$.

Now to get $f^{-1}(F)\in\M_{\a+1}(X\ti Y)$ it suffices to verify
that $f^{-1}(F)=\bigcap_{m=1}^\infty\bigcup_{n\ge m}A_{n,m}$.

To verify the inclusion $f^{-1}(F)\subset
\bigcap_{m\in\IN}\bigcup_{n\ge m}A_{n,m}$, fix any $m\in\IN$ and
$(x,y)\in X\ti Y$ with $f(x,y)\in F$. Let $U_n\in\U_n$, $n\in\IN$,
be a sequence of open sets containing the point $x$. Since
$s_n(U_n)\to x$ and $f$ is continuous with respect to the first
variable, we may find $k\ge m$ with $f(s_k(U_k),y)=f_{U_k}(y)\in
W_m$. Thus $(x,y)\in U_k\ti f_{U_k}^{-1}(W_m)\subset
A_{k,m}\subset \bigcup_{n\ge m}A_{n,m}$ for every $m$.

Next, assume that $(x,y)\notin f^{-1}(F)$. Then $f(x,y)\notin F$
and there is $k\in\IN$ with $\ol{W}_k\not\ni f(x,y)$. By the
choice of the covers $\U_n$ and the continuity of $f$ with respect
to the first variable, there is $m\ge k$ such that
$f(s_n(U_n),y)\notin \ol{W}_k$ for any $n\ge m$ and $x\in
U_n\in\U_n$. Since $W_m\subset W_k$, we get $(x,y)\notin
\bigcup_{n\ge m}A_{n,m}$.\qed
\enddemo

\remark{4.2. Remark} The inclusion $CC(X\ti Y,Z)\subset H_1(X\ti
Y,Z)$ is not true without restrictions on spaces $X,Y,Z$. If
$X=Y=\{\infty\}\cup\Gamma$ is the one-point compactification of un
uncountable discrete space $\Gamma$, then the function $f:X\ti
Y\to\IR$ defined by $$ f(x,y)=\cases 1,&\text{if
$x=y\in\Gamma$};\\ 0,&\text{otherwise}
\endcases
$$ is separately continuous but does not belong to the class
$H_1(X\ti Y,\IR)$.
\endremark

\heading 5. A generalization of the Rudin theorem
\endheading

In \cite{Ru} W.Rudin proved that every separately continuous
function $f:X\ti Y\to Z$ defined on the product of a metrizable
space $X$ and a topological space $Y$ and acting into a locally
convex topological vector space $Z$ belongs to the first Baire
class $B_1(X\ti Y,Z)$. In \cite{MMMS} it was proved that the
metrizability of $X$ Rudin Theorem can be replaced by the
paracompactness and the $\s$-metrizability; moreover, if the space
$X$ is finite-dimensional then the local convexity of $Z$ is
superfluous, see \cite{KM}.

In this section we shall prove that the Rudin Theorem is still
valid if $X$ is replaced by any metrically quarter-stratifiable
space and $Z$ by a locally convex equiconnected space $Z$.

We remind that an {\it equiconnected space} is a pair
$(Z,\lambda)$ consisting of a topological space $Z$ and a
continuous map $\lambda:Z\ti Z\ti[0,1]\to Z$ such that
$\lambda(x,y,0)=x$, $\lambda(x,y,1)=y$, and $\lambda(x,x,t)=x$ for
every $x,y\in X$, $t\in[0,1]$. For a subset $A\subset Z$ of an
equiconnected space $(Z,\lambda)$ let $\lambda^0(A)=A$ and
$\lambda^{n}(A)=\lambda(\lambda^{n-1}(A)\times A\ti[0,1])$ for
$n\ge 1$. Let also
$\lambda^\infty(A)=\bigcup_{n\in\IN}\lambda^n(A)$. An
equiconnected space $(Z,\lambda)$ is called {\it locally convex}
if for every point $z\in Z$ and a neighborhood $O(z)\subset Z$ of
$z$ there is a neighborhood $U\subset Z$ of $z$ with
$\lambda^\infty(U)\subset O(z)$.

Equiconnected spaces are tightly connected with absolute
extensors. We remind that a space $Z$ is an {\it absolute
extensors for a class $\C$} of topological spaces if every
continuous map $f:B\to Z$ from a closed subspace $B$ of a space
$C\in\C$ admits a continuous extension $\bar f:C\to Z$ over all
$C$. It is known that each equiconnected space is an absolute
extensor for strongly countable-dimensional stratifiable spaces,
while each locally convex equiconnected space is an absolute
extensor for the class of stratifiable spaces. Moreover, a
stratifiable space $Z$ is an absolute extensor for stratifiable
spaces if and only if $Z$ is admits an equiconnecting function
$\lambda$ turning $Z$ into a locally convex equiconnected spaces,
see \cite{Bo}, \cite{Ca}, \cite{BGS}.

Obvious examples of (locally convex) equiconnected spaces are
convex subsets of (locally convex) linear topological spaces and
their retracts. Also every contractible topological group $G$ is
an equiconnected space: an equiconnecting function $\lambda$ on
$G$ can be defined by $\lambda(x,y,t)=h_t(xy^{-1})\cdot
h_t(e)^{-1}\cdot y$, where $e$ is the neutral element of $G$ and
$\{h_t:G\to G\}_{t\in[0,1]}$ is a contraction of $G$ with
$h_0=\id$ and $h_1(G)=\{e\}$. Let us recall that a topological
space $X$ is {\it contractible} if it admits a homotopy
$h:X\ti[0,1]\to X$ such that $h(x,0)=x$ and $h(x,1)=*$ for all
$x\in X$ and some fixed point $*\in X$.

Given topological spaces $X$ and $Y$ denote by $C_p(X,Y)$ the
subspace of the Tychonov product $Y^X$, consisting of all
continuous functions from $X$ to $Y$. Let $B_0(X,Y)=C_p(X,Y)$ and
by transfinite induction for every ordinal $\a>0$ define the Baire
class $B_\a(X,Y)$ to be the sequential closure of the set
$\bigcup_{\xi<\a}B_\xi(X,Y)$ in $Y^X$.

Now let us look at the Rudin theorem from the following point of
view. Actually, this theorem states that every continuous function
$f:X\to C_p(Y,Z)$ from a metrizable space $X$ is a pointwise limit
of ``jointly continuous'' functions. This observation leads to the
following question: {\sl Under which conditions every continuous
map $f:X\to Z$ is a pointwise limit of some ``nice'' functions,
and what should be understood under a ``nice'' function?}

In case of an equiconnected space $(Z,\lambda)$ under ``nice''
functions we shall understand so-called piece-linear functions
which are defined as follows. Let $PL(X,Z)$ be the smallest
non-empty subset of $C_p(X,Z)$ satisfying the conditions: \roster
\item for every point $z_0\in Z$ and functions $f\in PL(X,Z)$, $\a\in
C(X,[0,1])$ the function $g(x)=\lambda(f(x),z_0,\alpha(x))$
belongs to $PL(X,Z)$;
\item a function $f:X\to Z$ belongs to $PL(X,Z)$ if there is an open
cover $\U$ of $X$ such that for every $U\in\U$ there is a function
$g\in PL(X,Z)$ with $f|U=g|U$.
\endroster

Equivalently, the set $PL(X,Z)$ can be defined constructively as
the set of functions $f:X\to Z$ for which there is an open cover
$\U$ of $X$ such that for every $U\in\U$ there are points
$z_0,\dots,z_n\in Z$ and functions $\a_1,\dots,a_n\in C(X,[0,1])$
such that $f|U=f_n|U$, where $f_0\equiv z_0$ and
$f_{i+1}(x)=\lambda(f_i(x),z_i,\alpha_i(x))$ for $0\le i<n$.

By $PL_{(1)}(X,Z)$ we denote the sequential closure of $PL(X,Z)$
in $Z^X$.

We recall that a topological space $X$ is defined to be {\it
strongly countable-dimensional} if $X$ can be represented as the
countable union $\bigcup_{n=1}^\infty X_n$ of closed subspaces
with $\dim X_n<\infty$ for all $n$. According to \cite{En$_2$,
5.1.10} a paracompact space $X$ is strongly countable-dimensional
if and only if $X$ can be represented as the countable union
$X=\bigcup_{n=1}^\infty X_n$ of closed subspaces such that every
open cover $\U$ of $X$ has an open refinement $\V$ with the
property $\ord\V|X_n\le n$ for every $n\in\IN$ (that is, every
point $X\in X_n$ belongs to at most $n$ elements of the cover
$\V$).

\proclaim{5.1. Theorem} Let $X$ be a metrically
quarter-stratifiable $T_1$-space and $(Z,\lambda)$ be an
equiconnected space. If $X$ is paracompact and strongly
countable-dimensional or $Z$ is locally convex, then
$C_p(X,Z)\subset PL_{(1)}(X,Z)$.
\endproclaim

\demo{Proof} Fix any function $f\in C_p(X,Z)$. Using Theorem 2.2,
find a weaker paracompact topology $\tau_p$ on $X$, a sequence
$\{\U_n\}_{n\in\IN}$ of $\tau_p$-open covers of $X$ and a sequence
$\{s_n:\U_n\to X\}_{n\in\IN}$ of functions such that $(x\in
U_n\in\U_n,\;n\in\IN)\Ra (s_n(U_n)\to x)$. Since the space
$(X,\tau_p)$ is paracompact, we may assume that each cover $\U_n$
is locally finite and contains no empty set.

If the space $X$ is paracompact and strongly
countable-dimensional, then according to the above-mentioned
Theorem 5.1.10 of \cite{En$_2$}, we may additionally assume that
$X$ is represented as the union $X=\bigcup_{m=1}^\infty X_m$ of
closed subsets such that $\ord\V_n|X_m\le m$ for every
$n,m\in\IN$.

For every $n\in\IN$ take any partition of unity
$\{\alpha_{n,U}:X\to[0,1]\}_{U\in\U_n}$ subordinated to $\U_n$,
that is $\a_{n,U}^{-1}(0,1]\subset U$ for each $U\in\U_n$ and
$\sum_{U\in\U_n}\a_{n,U}\equiv 1$.

Let $\xi_n:X\to N(\U_n)$,
$\xi_n:x\mapsto\sum_{U\in\U_n}\alpha_{n,U}(x)\cdot U$, be the
canonical map into the nerve of the cover $\U_n$, see \cite{FF,
VI.\S3}. Next, we construct a map $\eta_n:N(\U_n)\to Z$ as
follows. Denote by $N^{(k)}(\U_n)$ the $k$-skeleton of $N(\U_n)$,
where $k\ge 0$. The map $\eta_n$ will be defined by induction. Let
$\eta_n(U)=f(s_n(U))$ for every $U\in N^{(0)}(\U_n)=\U_n$. Suppose
that $\eta_n$ is already defined on the $k$-skeleton
$N^{(k)}(\U_n)$ of $N(\U_n)$. We shall extend $\eta_n$ onto
$N^{(k+1)}(\U_n)$. Take any point $x\in N^{(k+1)}(\U_n)\setminus
N^k(\U_n)$ and find a $k$-dimensional simplex $\s\ni x$. Fix any
vertex $v$ of $\s$. The point $x$ can be uniquely written as
$x=tv+(1-t)y$, where $t\in[0,1]$ and $y\in\sigma\cap
N^{(k)}(\U_n)$. Let $\eta_n(x)=\lambda(\eta_n(y), \eta_n(v),t)$
and observe that the so-extended map $\eta_n:N^{(k+1)}(\U_n)\to Z$
is continuous.

The inductive construction yields a continuous map
$\eta_n:N(\U_n)\to Z$ which has the following property:
$\eta_n(\s)\subset\lambda^\infty(\eta_n(\s^{(0)})=\lambda^\infty\big(\bigcup_{U\in\s^{(0)}}f(s_n(U))\big)$
for every simplex $\s$ of $N(\U_n)$. Observe also that
$\eta_n\circ\xi_n\in PL(X,Z)$.

We claim that $\eta_n\circ\xi_n\to f$, provided $X$ is strongly
countable-dimensional or $Z$ is locally convex. To show this, fix
any point $x\in X$ and a neighborhood $O(f(x))\subset Z$ of
$f(x)$.

First, we consider the case when the space $X$ is paracompact and
strongly countable-dimensional. In this case $x\in X_m$ for some
$m\in\IN$ and $\ord\U_n|X_m\le m$ for each $n$, which implies
$\xi_n(x)\in N^{(m)}(\U_n)$ for every $n$. Using the continuity of
the map $\lambda$, find a neighborhood $O_1\subset Z$ of $f(x)$
such that $\lambda^m(O_1)\subset O(f(x))$. By the choice of the
sequences $\{\U_n\}$ and $\{s_n\}$ there is $n_0\in\IN$ such that
$s_n(U_n)\in f^{-1}(O_1)$ for every $n\ge n_0$ and $U_n\in\U_n$
with $U_n\ni x$. It follows from the construction of the map
$\eta_n$ that $\eta_n\circ\xi_n(x)\in\lambda^m(O_1)\subset
O(f(x))$ for every $n\ge n_0$. Thus $\eta_n\circ\xi_n\to f$.

Now suppose that the equiconnected space $(Z,\lambda)$ is locally
convex. Then we may find a neighborhood $O_1\subset Z$ of $f(x)$
with $\lambda^\infty(O_1)\subset O(f(x))$ and a number $n_0\in\IN$
such that $s_n(U_n)\in f^{-1}(O_1)$ for every $n\ge n_0$ and $x\in
U_n\in\U_n$. It follows from the construction of the map $\eta_n$
that $\eta_n\circ \xi_n(x)\subset \lambda^\infty(O_1)\subset
O(f(x))$ for every $n\ge n_0$. Thus in both cases the function
sequence $\{\eta_n\circ\xi_n\}_{n=1}^\infty\subset PL(X,Z)$ tends
to $f$, which proves that $f\in PL_{(1)}(X,Z)$.\qed
\enddemo

Now we use Theorem 5.1 to generalize the Rudin Theorem as well as
the results of \cite{MMMS}, \cite{MS}, \cite{KMM}, \cite{KM}. By
$CC(X\ti Y,Z)$ the set of separately continuous functions $X\ti
Y\to Z$ is denoted.

\proclaim{5.2. Corollary} Let $X$ be a metrically
quarter-stratifiable space, $Y$ be a topological space, and
$(Z,\lambda)$ be an equiconnected space. If $X$ is paracompact and
strongly countable-dimensional or $Z$ is locally convex, then
$CC(X\ti Y,Z)\subset B_1(X\ti Y,Z)$.
\endproclaim

\demo{Proof} In obvious way the equiconnecting function $\lambda$
induces an equiconnecting function on $C_p(Y,Z)$. Moreover, the
equiconnected space $C_p(Y,Z)$ is locally convex if so is $Z$.
Hence, we may apply Theorem 5.1 to conclude that
$C_p(X,C_p(Y,Z))\subset PL_{(1)}(X,C_p(Y,Z))$. Observe that the
space $C_p(X,C_p(Y,Z))$ may be identified with $CC(X\ti Y,Z)$,
while every element of $PL(X,C_p(Y,Z))$ is jointly continuous as a
function $X\ti Y\to Z$. This yields that $CC(X\ti Y,Z)\subset
B_1(X\ti Y,Z)$.\qed
\enddemo

For an ordinal $\a$ and topological spaces $X,Y,Z$ denote by
$CB_\a(X\ti Y,Z)$ the set of functions $f:X\ti Y\to Z$ such that
for every $x_0\in X$ and $y_0\in Y$ we have $f(\cdot,y_0)\in
C(X,Z)$ and $f(x_0,\cdot)\in B_\a(Y,Z)$. Thus $CB_0(X\ti
Y,Z)=CC(X\ti Y,Z)$. Generalizing the Rudin Theorem,
V.K.Maslyuchenko et al. \cite{MMMS}
 proved that $CB_\a(X\ti Y,Z)\subset
B_{\a+1}(X\ti Y,Z)$ for any countable ordinal $\a$, any metrizable
space $X$, any topological space $Y$ and any locally convex space
$Z$, see \cite{MMMS}. Below we have a further generalization.

\proclaim{5.3. Theorem} Let $X$ be a metrically
quarter-stratifiable space, $Y$ be a topological space and $Z$ be
a contractible space. Then $CB_\a(X\ti Y,Z)\subset B_{\a+1}(X\ti
Y,Z)$ for every countable ordinal $\a>0$.
\endproclaim

\demo{Proof} Let $\a\ge 1$ be a countable ordinal and $f\in
CB_\a(X\ti Y,Z)$. Fix a non-decreasing sequence $(\a_n)$ of
ordinals with $\a=\sup_n\a_n+1$ and let $h:Z\ti[0,1]\to Z$ be a
contraction of $Z$ with $h(z,0)=z$ and $h(z,1)=*$ for all $z\in Z$
and some fixed point $*\in Z$.

Using Theorem 2.2, find a weaker metrizable topology $\tau_m$ on
$X$,
 a sequence $\{\U_n\}_{n\in\IN}$ of $\tau_m$-open cover of
$X$ and a sequence $\{s_n:\U_n\to X\}_{n\in\IN}$ of maps such that
$s_n(U_n)\to x$ for every $x\in U_n\in\U_n$, $n\in\IN$. Using the
paracompactness of $(X,\tau_m)$, for every $n\in\IN$ find a
partition of unity $\{\a_{n,U}:(X,\tau_m)\to[0,1]\}_{U\in\U_n}$
subordinated to $\U_n$, that is, $\alpha_{n,U}^{-1}(0,1]\subset U$
for $U\in\U_n$. Then $\{\a_{n,U}^{-1}(0,1]\}_{U\in\U_n}$ is a
locally finite cover of $X$ by functionally open sets.

Fix any well-ordering $\le$ on the set $\U_n$. For every
$U\in\U_n$ and $m\in\IN$ consider the closed set
$F_{m,U}=\a_{n,U}^{-1}[1/m,1]\setminus\bigcup_{V<U}\a_{n,V}^{-1}(0,1]$.
Clearly, $\{F_{m,U}\}_{U\in\U_n}$ is a discrete collection of
closed subsets of $(X,\tau_m)$ for every $m\in\IN$. Using the
collective normality of $(X,\tau_m)$, for every
$(m,U)\in\IN\ti\U_n$ we may find a continuous function
$\beta_{m,U}:(X,\tau_m)\to [0,1]$ such that
$\beta_{m,U}^{-1}(1)\supset F_{m,U}$ and
$\{\beta_{m,U}^{-1}(0,1]\}_{U\in\U_n}$ is a discrete collection in
$(X,\tau_m)$ for every $m\in\IN$.

For every $U\in\U_n$ find a function sequence
$\{f_{m,U}\}_{m\in\IN}\subset\bigcup_{\xi<\a}B_\xi(Y,Z)$
with\break $\lim_{m\to\infty}f_{m,U}(y)=f(s_n(U),y)$ for every
$y\in Y$. Without loss of generality, we may assume that $f_m\in
B_{\a_m}(Y,Z)$ for every $m$. Then the function $g_{m,n}:X\ti Y\to
Z$ defined by the formula $$ g_{n,m}(x,y)=\cases h(f_{m,U}(y),
1-\beta_{m,U}(x)), &\text{if $\beta_{m,U}(x)>0$ for some
$U\in\U_n$};\\ *,&\text{otherwise}
\endcases
$$ belongs to the class $B_{\a_m}(X\ti Y,Z)$ for every
$m,n\in\IN$. Let us show that the limit
$\lim_{m\to\infty}g_{m,n}(x,y)$ exists for every $(x,y)\in X\ti Y$
and $n\in\IN$. Indeed, given $(x,y)\in (X,Y)$ and $n\in\IN$ let
$U_n=\min\{U\in\U_n:\a_{n,U}(x)>0\}$. Then $x\in F_{m,U_n}$ for
all sufficiently large $m$ and
$g_{m,n}(x)=h(f_{m,U_n}(y),1-\beta_{m,U_n}(x))=f_{m,U_n}(y)\to
f(s_n(U_n),y)$ as $m\to\infty$. Thus
$g_n(x,y)=\lim_{m\to\infty}g_{m,n}(x,y)=f(s_n(U_n),y)$ exists and
$g_n\in B_\a(X\ti Y,Z)$. Next, since $s_n(U_n)\to x$, we get
$f(s_n(U_n),y)\to f(x,y)$ which yields that $f=\lim_{n\to\infty}
g_n$ and $f\in B_{\a+1}(X\ti Y,Z)$.\qed
\enddemo

The following corollary generalizes a result from \cite{MMMS}.

\proclaim{5.4. Corollary} Let $X_1,\dots,X_n$ be metrically
quarter-stratifiable spaces, $Y$ be a topological space, and
$(Z,\lambda)$ be an equiconnected space. If $X_1$ is paracompact
and strongly countable-dimensional or $Z$ is locally convex, then
every separately continuous function $f:X_1\times\dots\times
X_n\times Y\to Z$ belongs to the $n$-th Baire class.
\endproclaim

We do not know if the conditions on $X_1$ or $Z$ are essential.

\proclaim{5.5. Question} Does every separately continuous function
$f:X\ti Y\to Z$ defined on the product of metrizable compacta and
acting into a linear metric space belong to the first Baire class?
\endproclaim

\definition{5.6. Remark}
Actually, according to \cite{MMMS}, the Rudin Theorem holds in a
more strong form: $\bar CC(X\times Y,Z)\subset B_1(X\ti Y,Z)$ for
any metrizable space $X$, topological space $Y$, and a locally
convex space $Z$, where $\bar CC(X\ti Y,Z)$ stands for the set of
functions for which there is a dense subset $D\subset X$ such that
for every $x_0\in D$ and every $y_0\in Y$ the functions
$f(x_0,\cdot)$ and $f(\cdot,y_0)$ are continuous. It is
interesting to remark that the metrizability of $X$ is essential
in this stronger form of the Rudin Theorem and can not be replaced
by stratifiability; namely $\bar CC(A\times[0,1],\IR)\not\subset
B_1(A\times[0,1],\IR)$, where $A$ is the Arens fan.
\enddefinition

We recall the definition of the {\it Arens fan} $A$, see
\cite{En$_1$, 1.6.19}.  Consider the following points of the real
line: $a_0=0$, $a_n=\frac1n$ and $a_{nm}=\frac1n+\frac1m$, where
$1\le n^2\le m$. Let $A_k=\{a_n:n\ge 0\}\cup\{a_{nm}:k\le n\le
n^2\le m\}$ for $k\in\IN$. On the union $A=\bigcup_{k\in\IN}A_k$
consider the strongest topology inducing the original topology on
each compactum $A_k$. Clearly, the so-defined space $A$ is
countable and stratifiable. The space $A$ is known as an example
of a sequential space which is not a Fr\'echet-Urysohn space.

\proclaim{5.7. Example} $\bar CC(A\times [0,1],\IR)\not\subset
B_1(A\times [0,1],\IR)$ for the Arens fan $A$.
\endproclaim

\demo{Proof} We shall construct a function $f\in\bar
CC(A\ti[0,1],\IR)$ which is not of the first Baire class. Take any
function $f_0\in B_2([0,1],\IR)\setminus B_1([0,1],\IR)$. Write
$f_0=\lim_{n\to\infty}f_n$ where $\{f_n\}_{n\in\IN}\subset
B_1([0,1],\IR)$. In its turn, represent each $f_n$ as a pointwise
limit $f_n=\lim_{m\to\infty}f_{nm}$ of continuous functions. Now
consider the map $f:A\times [0,1]\to \IR$: $$ f(a,y)=\cases
f_0(y),&\text{ if $a=a_0$};\\ f_n(y),&\text{ if $a=a_n$};\\
f_{nm}(y),&\text{ if $a=a_{nm}$.}
\endcases
$$ It is easy to see that $f\in\bar CC(A\ti[0,1],\IR)$ but
$f\notin B_1(A\ti[0,1],\IR)$. \qed
\enddemo

\heading 6. Rudin spaces
\endheading

In light of the mentioned generalizations of the Rudin Theorem it
is natural to introduce the following

\definition{6.1. Definition} A topological space $X$ is defined to be {\it
Rudin} if for arbitrary topological space $Y$ every separately
continuous function $f:X\ti Y\to\IR$ belongs to the first Baire
class.
\enddefinition

In the following theorem we collect all fact concerning Rudin
spaces we know at the moment. Let us recall that a subspace $Y$ of
a topological space $X$ is {\it $t$-embedded}, if there is a
continuous extender $E:C_p(Y)\to C_p(X)$, that is a map $E$ such
that $E(f)|Y=f$ for every $f\in C_p(Y)$, see \cite{Ar$_1$,
III.\S2.1}. Here $C_p(X)=C_p(X,\IR)$.

\proclaim{6.2. Theorem} \roster
\item A space $X$ is Rudin if and only if the calculation map $c_X:X\ti
C_p(X)\to \IR$, $c_X:(x,f)\mapsto f(x)$, is of the first Baire
class.
\item A space $X$ is Rudin if the space $C_p(X)$ is Rudin.
\item Every metrically quarter-stratifiable space is Rudin.
\item If a space $X$ is Rudin, then each $t$-embedded subspace of $X$ is
Rudin.
\item A space $X$ is Rudin if $X=A\cup B$, where $A,B$ are  Rudin
closed subspaces of $X$ and $A$ is a functionally closed retract
of $X$.
\item If a Tychonov space $X$ is Rudin, then $d(X)\le l(X)$. Consequently,
each compact Rudin space is separable.
\endroster
\endproclaim

\demo{Proof} 1) The ``only if'' part of the first statement is
trivial. To prove the ``if'' part, assume that the calculation map
$c:X\ti C_p(X)\to \IR$ is of the first Baire class. Observe that
every separately continuous map $f:X\ti Y\to\IR$ can be seen as a
continuous map $F:Y\to C_p(X)$ such that $f(x,y)=c(x,F(y))$ for
$(x,y)\in X\ti Y$. Now it is clear that $f$ is of the first Baire
class. \vskip5pt

2) The second statement follows immediately from the first one.
\vskip5pt

3) According to Corollary 4.2, each metrically
quarter-stratifiable space is Rudin. \vskip5pt

4) Suppose $X$ is a Rudin space and $Y$ a $t$-embedded subspace of
$X$. Let $E:C_p(Y)\to C_p(X)$ be the corresponding extender. Let
$c_X:X\ti C_p(X)\to \IR$, $c_Y:Y\times C_p(Y)\to\IR$ be the
calculation functions for spaces $X$ and $Y$, respectively. Taking
into account that the function $c_X$ is of the first Baire class
and $c_Y(y,f)=c_X(y,E(f))$ for every $(y,f)\in Y\times C_p(Y)$ we
conclude that $c_Y$ is a function of the first Baire class either.
\vskip5pt

5) Suppose $X=A\cup B$, where $A,B$ are Rudin closed subspaces of
$X$ and $A$ is a functionally closed retract in $X$. Let $r:X\to
A$ be a retraction and $\psi:X\to[0,1]$ a map with
$\psi^{-1}(0)=A$. Since the spaces $A,B$ are Rudin, there are
sequences of maps $\a_n:A\times C_p(A)\to\IR$ and $\beta_n:B\ti
C_p(B)\to\IR$, $n\in\IN$, tending to the calculation functions
$c_A$ and $c_B$ of the spaces $A,B$, respectively. For every
$n\in\IN$ consider the function $\psi_n=\min\{n\psi,1\}$ and
define the map $c_n:X\ti C_p(X)\to\IR$ by the formula $$
c_n(x,f)=\cases \alpha_n(x,f|A),& \text{if $x\in A$};\\
(1-\psi_n(x))\alpha_n(r(x),f|A)+\psi_n(x)\beta_n(x,f|B),&\text{if
$x\in B$}.
\endcases
$$ It is easy to see that the maps $c_n$ are well-defined and
continuous, and the sequence $\{c_n\}$ tends to the calculation
function of $X$. \vskip5pt

6) Suppose $X$ is a Rudin space. Let $\{c_n\}_{n\in\IN}\subset
C(X\ti C_p(X),\IR)$ be a sequence of continuous functions tending
to the calculation function $c:X\ti C_p(X)\to\IR$. Then
$c^{-1}(0)=\bigcap_{n\in\IN}\bigcup_{m\ge n}c_m^{-1}(-1/n,1/m)$ is
a $G_\delta$-set in $X\ti C_p(X)$. Let $\{O_n\}_{n\in\IN}$ be a
sequence of open sets in $X\ti C_p(X)$ with
$\bigcap_{n=1}O_n=c^{-1}(0)$. By the continuity of the functions
$c_n$, for every $x\in X$ and $n\in\IN$ we may find a neighborhood
$U_n(x)$ and a finite subset $F_n(x)\subset X$ such that
$U_n(x)\ti F_n(x)^\perp\subset O_n$, where $A^\perp=\{f\in
C_p(X):f|A\equiv0\}$ for a subset $A\subset X$. For every $n\in
\IN$ choose a subset $X_n\subset X$ of size $|X_n|\le l(X)$ with
$\bigcup_{x\in X_n}U_n(x)=X$. We claim that $D=\cup\{F_n(x):x\in
X_n,\; n\in\IN\}$ is a dense set in $X$. Otherwise, we would find
a function $f\in C_p(X)$ with $f|D\equiv 0$ and $f(x)\ne 0$ for
some $x\in X$. It is easy to see that $(x,f)\in
\bigcap_{n\in\IN}O_n$ but $(x,f)\notin c^{-1}(0)$, a
contradiction, which shows that $D$ is dense in $X$ and
$d(X)\le|D|\le l(X)$.\qed
\enddemo

There are many open question about the structure of Rudin spaces.

\definition{6.3. Question} Is there a Rudin space which is not
quarter-stratifiable? Is a (compact) space $X$ Rudin, if
$\{(x,f):f(x)=0\}$ is a $G_\delta$-set in $X\ti C_p(X)$? Is every
(fragmentable, scattered or Rosenthal) compact Rudin space
metrizable?
\enddefinition

Let us remark that according to \cite{Ar$_2$, II.6.1} each Corson
Rudin compactum, being separable, is metrizable. Concerning the
last question let us remark that there are non-metrizable compact
spaces $X$ such that $CC(X\times Y,\IR)\subset B_1(X\times Y,\IR)$
for every compact space $Y$. The following result belongs to
G.Vera \cite{Ve}.

\proclaim{6.4. Theorem} A compact space $X$ has countable Souslin
number if and only if for every compact space $Y$ every separately
continuous function $f:X\times Y\to\IR$ is of the first Baire
class.
\endproclaim

\demo{Proof} Suppose $X$ has a countable Souslin number and
$f:X\times Y\to \IR$ is a separately continuous function which can
be seen as a continuous map $F:Y\to C_p(X)$. Then the image
$F(Y)\subset C_p(X)$ is a compact subset of the space $C_p(X)$
over a compactum with countable Souslin number. According to
Theorem 14 \cite{Ar$_1$, II.\S2} all such compact subsets are
metrizable. By the Rudin Theorem, the restriction of the
separately continuous calculation function $c:X\times
C_p(X)\to\IR$ onto $X\ti F(Y)$ is of the first Baire class, which
implies that the map $f(x,y)=c(x,F(y))$ is of the first Baire
class either.

Now assume that the Souslin number of a compactum $X$ is
uncountable. Then there exists an uncountable family $\U$ of
pairwise disjoint nonempty open sets in $X$. For every $U\in\U$
fix a continuous function $f_U:X\to[0,1]$ such that
$f_U|X\setminus U\equiv 0$ and $\max f_U=1$. Denote by $0$ the
origin of the vector space $C_p(X)$. It is easy to see that the
set $Y=\{0\}\cup\{f_U:U\in\U\}$ in $C_p(X)$ is compact and the
calculation function $c:X\ti Y\to\IR$, $c:(x,f)\mapsto f(x)$, is
separately continuous. Yet, since $c^{-1}(1/2,1]$ is not an
$F_\sigma$-set in $X\ti Y$, the function $c$ is not of the first
Baire class. In fact, one may show that
$c\notin\bigcup_{\a<\w_1}B_\a(X\ti Y,Z)$, i.e., $c$ is not Baire
measurable. \qed
\enddemo

\definition{6.5. Remark}
It follows from Theorems 6.2 and 6.4 that any unseparable compact
space $X$ with countable Souslin number is not Rudin but for every
compact space $Y$ every separately continuous function $f:X\ti
Y\to\IR$ belongs to the first Baire class.
\enddefinition

Finally let us prove one more result, connected with Baire
classification of separately continuous functions.

\proclaim{6.6. Theorem} Let $X,Y,Z$ be compact Hausdorff spaces.
If $X$ is separable and $Y$ has countable Souslin number, then
every separately continuous function $f:X\ti Y\ti Z\to\IR$ belongs
to the second Baire class.
\endproclaim

\demo{Proof} Fix any separately continuous function $f:X\ti Y\ti
Z\to\IR$, which can be seen as a continuous map $F:X\to CC(Y\ti
Z,\IR)$ into the space of separately continuous functions endowed
with the product topology. According to a recent result of S.Gulko
and G.Sokolov \cite{GS} if the space $Y$ has countable Souslin
number, then every compactum in $CC(Y\times Z,\IR)$ is Corson. Now
if $X$ is separable and $Y$ has countable Souslin number, then
$F(X)\subset CC(Y\ti Z,\IR)$ is a separable Corson compactum.
Since separable Corson compacta are metrizable, see \cite{Ar$_2$,
II.6.1}, we conclude that the compactum $F(X)$ is metrizable. By
\cite{MMMS} (see also Corollary 5.4), the calculation map
$c:F(X)\ti Y\ti X\to\IR$, $c:(g,y,z)\mapsto g(y,z)$, being
separately continuous, belongs to the second Baire class. Then the
function $f$ belongs to the second Baire class because
$f(x,y,z)=c(F(x),y,z)$.\qed
\enddemo

\definition{6.6. Remark} The class of Rudin spaces includes the class of
Lebesgue spaces introduced by O.Sobchuk in \cite{So}.
\enddefinition

\enddocument